\documentclass[11pt,reqno]{amsart}

\usepackage[dvipsnames]{xcolor}
\usepackage{amssymb}
\usepackage{mathrsfs}
\usepackage{amsmath}
\usepackage{amsthm}
\usepackage{url}
\usepackage{stackrel} 
\usepackage{courier}
\usepackage[all,cmtip]{xy} 
\usepackage{multicol} 
\usepackage{verbatim} %  useful for \begin{comment} \end{comment}  
\usepackage{mathtools}  
\usepackage{hyperref}
\usepackage{multicol}  
\usepackage{enumerate} 
\usepackage{bookmark} 
\usepackage{ytableau}
\usepackage{microtype}
 
\catcode`~=11 % make LaTeX treat tilde (~) like a normal character
\newcommand{\urltilde}{\kern -.15em\lower .7ex\hbox{~}\kern .04em}  
\catcode`~=13 % revert back to treating tilde (~) as an active character

\newcommand{\CC}{\mathbb C} 

\newcommand{\Lie}{\mathop{\mathrm{Lie}}\nolimits}

\newcommand{\lie}{\mathop{\mathrm{Lie}}\nolimits}

\newcommand{\op}{\mathop{\mathrm{op}}\nolimits}

\newcommand{\Proj}{\mathop{\mathrm{Proj}}\nolimits}

\newcommand{\codim}{\mathop{\mathrm{codim}}\nolimits}

\newcommand{\I}{\mathop{\mathrm{I}_n}\nolimits}

\newcommand{\spec}{\mathop{\mathrm{Spec}}\nolimits}  

\newcommand{\tr}{\mathop{\mathrm{tr}}\nolimits}

\newcommand{\Tr}{\mathop{\mathrm{tr}}\nolimits}

\newcommand{\Rep}{\mathop{\mathrm{Rep}}\nolimits} 
\newcommand{\Hom}{\mathop{\mathrm{Hom}}\nolimits}

\newcommand{\GL}{\mathop{\mathrm{GL}}\nolimits} 
 
\newcommand{\rss}{\mathop{\mathrm{rss}}\nolimits} 
\newcommand{\End}{\mathop{\mathrm{End}}\nolimits}

\makeatletter
\newif\if@borderstar
   \def\bordermatrix{\@ifnextchar*{%
       \@borderstartrue\@bordermatrix@i}{\@borderstarfalse\@bordermatrix@i*}%
   }
   \def\@bordermatrix@i*{\@ifnextchar[{\@bordermatrix@ii}{\@bordermatrix@ii[()]}}
   \def\@bordermatrix@ii[#1]#2{%
   \begingroup
     \m@th\@tempdima8.75\p@\setbox\z@\vbox{%
       \def\cr{\crcr\noalign{\kern 2\p@\global\let\cr\endline }}%
       \ialign {$##$\hfil\kern 2\p@\kern\@tempdima & \thinspace %
       \hfil $##$\hfil && \quad\hfil $##$\hfil\crcr\omit\strut %
       \hfil\crcr\noalign{\kern -\baselineskip}#2\crcr\omit %
       \strut\cr}}%
     \setbox\tw@\vbox{\unvcopy\z@\global\setbox\@ne\lastbox}%
     \setbox\tw@\hbox{\unhbox\@ne\unskip\global\setbox\@ne\lastbox}%
     \setbox\tw@\hbox{%
       $\kern\wd\@ne\kern -\@tempdima\left\@firstoftwo#1%
         \if@borderstar\kern2pt\else\kern -\wd\@ne\fi%
       \global\setbox\@ne\vbox{\box\@ne\if@borderstar\else\kern 2\p@\fi}%
       \vcenter{\if@borderstar\else\kern -\ht\@ne\fi%
         \unvbox\z@\kern-\if@borderstar2\fi\baselineskip}%
         \if@borderstar\kern-2\@tempdima\kern2\p@\else\,\fi\right\@secondoftwo#1 $%
     }\null \;\vbox{\kern\ht\@ne\box\tw@}%
   \endgroup
   }
\makeatother

\newtheorem{theorem}{Theorem}[section]

\newtheorem{problem}[theorem]{Problem}

\newtheorem{assumption}[theorem]{Standing Assumption}

\allowdisplaybreaks 

\keywords{Nakajima quiver variety, geometric invariant theory, nonreductive group, filtered ADHM data, Hilbert scheme, rational Cherednik algebra}    
\subjclass[2000]{Primary:  16G20,  	20G05. Secondary:  	20G20, 14C05}

\begin{document}

\title[Affine and GIT quotients of the extended Grothendieck--Springer resolution]{Suggestions to study affine and GIT quotients of the extended Grothendieck--Springer resolution}  
\author{Mee Seong Im} 
\address{Department of Mathematical Sciences, United States Military Academy, West Point, NY 10996 USA}
\email{meeseongim@gmail.com}  
\date{\today}

\begin{abstract}   
We define filtered ADHM data and connect a notion of filtered quiver representations to Grothendieck--Springer resolutions. 
We also provide current developments and give a list of research problems to further study filtered ADHM equation. 
\end{abstract}   

\maketitle

\setcounter{tocdepth}{1}

\section{Introduction}
\label{section:introduction}  
Springer resolution and the Grothendieck--Springer resolution are fundamental and important in representation theory. 
In type A, the Springer correspondence gives a bijection between irreducible symmetric group representations and unipotent conjugacy classes of the general linear group.  That is, given a unipotent conjugacy class and a fixed element $u$ in the conjugacy class, the corresponding irreducible representation of $S_n$ is the cohomology group $H^{2\dim \mathcal{B}_u}(\mathcal{B}_u,\mathbb{Q})$, where $\mathcal{B}_u$ is the set of Borel subgroups of $GL_n(\mathbb{C})$ in the Springer resolution of the unipotent group containing $u$ (cf. \cite{MR618892,MR2838836,MR1649626}). 

 Embedded in the Grothendieck--Springer resolution is the Springer resolution (cf. \cite{MR2838836,Im-rss-locus-GS-resolution}), and 
the enhanced Grothendieck--Springer resolution could be viewed as the Hamiltonian reduction of a certain parabolic moment map (cf.~\cite{Nevins-GSresolutions,Im-rss-locus-GS-resolution}), where the latter is  easier to visualize, manipulate and study, using a notion called filtered quiver representations 
(cf.~\S\ref{subsubsection:filtered-ADHM-data} and \S\ref{subsubsection:quiver-varieties}; also see~\cite[Prop.~3.2]{Nevins-GSresolutions}, \cite[Prop.~1.1]{Im-rss-locus-GS-resolution}).

Also appearing in representation theory, quiver varieties arise in many context and have deep connections to mathematical physics, string theory, cluster algebras, Kac--Moody algebras, to name a few. Lusztig's \cite{MR1775358,MR1088333,MR1623674} 
and Nakajima's \cite{MR1302318,MR1711344} are some of many foundational grounds to study quiver representations. 
In this paper, we generalize Nakajima quiver varieties (cf.~\S\ref{subsubsection:filtered-ADHM-data}), provide current development (as a generalization of classical results), and give a list of open problems (cf.~\S\ref{section:results-problems}). 
Current development includes using a new technique through filtered ADHM data, 
where we explicitly give in \S\ref{section:results-filtered-semi-invariants} 
(semi-)invariant polynomials for certain finite acyclic and cyclic quiver representations with a fixed filtration (see \cite{Im-two-pathways,Im-doctoral-thesis} for the strategies). 
We then focus on the Grothendieck--Springer resolution setting and describe the 
$B$ (and $P$-)orbits on the Hamiltonian reduction of the cotangent bundle of extended $\mathfrak{b}=\Lie(B)$ 
(and $\mathfrak{p}=\Lie(P)$) in \S\ref{section:Ham-red-Borel} and \S\ref{section:Ham-red-parabolic}, respectively.   
In \S\ref{section:results-problems}, 
we include constructing affine and GIT quotients using filtered ADHM equations (cf.~\S\ref{subsection:affine-GIT-open}), 
discovering connections to modified rational Cherednik algebras (cf.~\S\ref{subsection:rational-Cherednik-alg}) and corresponding $\mathcal{D}$-modules, 
and investigating variations of Hilbert schemes (cf.~\S\ref{subsection:Hilbert-schemes}). 
Some motivations to generalize quiver varieties come from the study of quiver Hecke (KLR-)algebras (cf.~\cite{MR2525917,MR2763732,Rouquier-2-Kac-Moody-algebras}), 
universal quiver Grassmannians and universal quiver flag varieties 
(cf.~\cite{MR2772068}),  
Lusztig's triangular decomposition of the upper half 
  of the universal enveloping algebra of a Kac--Moody algebra  
(cf.~\cite{MR1035415,MR1182165,MR1758244}), and 
 the Grothendieck--Springer resolution (cf.~\cite{MR2838836,Nevins-GSresolutions}).

The construction of Nakajima quiver varieties involves an action by reductive groups, but in our generalized setting, we shift our focus to parabolic group actions.

Studying parabolic group actions on vector spaces with a fixed filtration 
and connecting affine and GIT constructions to well-known varieties,    
revealing different ways to view the same geometric object and often providing a deeper insight to known varieties, are important in mathematics. 
First, 
consider a classical example: let $G=\GL_n(\mathbb{C})$ and $\mathfrak{g}=\mathfrak{gl}_n=\Lie(G)$.   
In the case of $G$-action on its Lie algebra $\mathfrak{g}$ by conjugation,    
the Jordan quiver in \S\ref{subsubsection:quiver-rep} is the natural quiver associated to this $G$-equivariant geometry.   
The orbit space $\mathfrak{g}/G$ is not a variety, while the affine quotient $\mathfrak{g}/\!\!/G$ consists of equivalence classes of semisimple, i.e., diagonalizable, matrices, 
so it is isomorphic to $\mathbb{C}^n$. Algebrically,  
$\mathfrak{g}/\!\!/G := \spec( \mathbb{C}[\mathfrak{g}]^{G})$,  
with the generators being coefficients of the characteristic polynomial of $n\times n$ matrices. 
Now, let $B\leq G$ be the standard Borel subgroup, i.e., the set of invertible upper triangular matrices, and let $\mathfrak{b}=\lie(B)$.  
The $B$-action on $\mathfrak{b}$ is naturally a quiver representation but with a notion of filtration modifying the representation space of the Jordan quiver,   
which one may think of $\text{F}^{\bullet}\Rep(Q,n)$ to be the set of all those   
 linear maps $W_r\in \End(\mathbb{C}^n)$   
 fixing the complete standard flag $\text{F}^{\bullet}$ in $\mathbb{C}^n$. 
 The subspace $\text{F}^{\bullet}\Rep(Q,n)$ is acted upon by those matrices in $G$ which preserve this flag, i.e, this is precisely $B$ (cf.~\S\ref{subsubsection:filtered-ADHM-data}).

 Throughout this paper, we will work over $\mathbb{C}$, and we will assume that 
 the set $\mathbb{N}$ of natural numbers includes $0$. 
 
\subsection*{Acknowledgment}   
The author would like to thank Thomas Nevins for invaluable discussions in representation theory during the initial stages of this paper. 
The author would also like to thank Darlayne Addabbo and Bolor Turmunkh for productive research meetings to discuss topics in geometric representation theory, and, in particular, 
geometric constructions of numerous important algebras and on Nakajima's work, in summer of 2014.  
M.S.I. was partially supported by NSA grant H98230-12-1-0216, by Campus Research Board at the University of Illinois at Urbana-Champaign, and by NSF grant DMS 08-38434.  

\section{Representation theory}
\label{section:rep-thy}

\subsection{Parabolic equivariant geometry}
\label{subsection:intro-equivariant-geometry}  
Studying $GL_n(\mathbb{C})$-orbits on the set $\mathfrak{gl}_n$ of $n\times n$ complex matrices by conjugation amounts to putting the 
matrices into Jordan canonical form, up to a permutation of its elementary blocks.  
In such area of mathematics, one hopes to construct interesting moduli spaces associated to a pair $(\mathfrak{X},\mathfrak{G})$, where $\mathfrak{X}$ is a space and $\mathfrak{G}$ is a group. 
%Such construction of moduli spaces will be discussed in Section~\ref{subsection:background-quivers-GIT}.   

Let $B$ be invertible upper triangular matrices in $G := \GL_n(\mathbb{C})$  
and let $\mathfrak{b}$ be the set of upper triangular matrices in $\mathfrak{g}=\mathfrak{gl}_n$. Letting $B$ act on $\mathfrak{b}$ by adjoint action, we can ask what are the $B$-orbits on $\mathfrak{b}$?   
More generally, how does one manage a general parabolic group action on a variety that has a certain notion of filtrations associated to it?

Returning to the example of the $B$-action on $\mathfrak{b}$, 
such pair $(\mathfrak{b},B)$ can be thought of as 
a variety with a natural notion of filtration.  
The set $\mathbb{C}^n \stackrel{a}{\rightarrow}\mathbb{C}^n$ of all linear homomorphisms can be identified with $\mathfrak{g}$,  
and changing domain and codomain basis vectors amounts to a certain action by $\GL_n(\mathbb{C})$. 
Now suppose we impose a complete filtration on both the domain and the codomain with respect to the standard basis vectors. That is, 
letting $\text{F}^{\bullet}: 0 \subseteq \mathbb{C}^1\subseteq \mathbb{C}^2\subseteq \ldots \subseteq \mathbb{C}^n$ be the complete standard filtration of vectors, where $\mathbb{C}^k$ is the subspace of $\mathbb{C}^n$ spanned by the first $k$ standard basis vectors, 
consider the set $\text{F}^{\bullet}\stackrel{r}{\rightarrow} \text{F}^{\bullet}$ of all linear maps preserving the filtration;  
such maps consist of $r$ such that $r|_{\mathbb{C}^k}:\mathbb{C}^k\rightarrow \mathbb{C}^k$ is linear for each $k$. 
Since $\text{F}^{\bullet}$ is a filtration with respect to standard basis vectors, $r$ has the form of an upper triangular matrix,  and changing domain and codomain basis vectors (which preserves the filtration) amounts to an action by $B$. 

\subsubsection{Quiver representations}
\label{subsubsection:quiver-rep}
We refer to  
\cite{Crawley-Boevey-rep-quivers, Ginzburg-Nakajima-quivers,Im-doctoral-thesis} 
for an extensive background on quivers and their representations.    
Quivers $Q$ are a directed graph with finite number of vertices and arrows.   
We will denote the set of vertices as $Q_0$ and the set of arrows as $Q_1$.  
 An example is the Jordan quiver: 
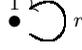
\[
  \xymatrix{   \stackrel{1}{\bullet} \ar@(rd,ru)_{r} 
	       }   
\] 
with one vertex $1$ and one arrow $r$; in this example, the head $hr$ of $r$ and the tail $tr$ of $r$   
end and begin at the same vertex.  We call the Jordan quiver cyclic since it has an oriented cycle, i.e., 
$r$ is called a cycle since the tail of $r$ equals the head of $r$.    
Another example is the $A_2$-quiver:    
\[   
 \xymatrix{  \stackrel{1}{\bullet} \ar[rr]^{a} & & \stackrel{2}{\bullet}   
}    
\]  
with two vertices $1$ and $2$ and one arrow $a$ whose tail of $a$ equals $1$ and the head of $a$ equals $2$. 
We call this quiver the $A_2$-quiver   
since the underlying graph has the structure of an $A_2$-Dynkin diagram. 
Furthermore, we call the $A_2$-quiver acyclic since it does not have any oriented cycles.

Let $\mathbf{v}=( v_1, \ldots, v_{Q_0} ) \in \mathbb{N}^{Q_0}$, 
a dimension vector. 
A representation $V$ of a quiver with dimension vector $\mathbf{v}$ 
assigns a finite dimensional vector space $V_i$ of dimension $v_i$ 
to each vertex and a linear map $V_a$  
to each arrow $a\in Q_1$.   
The representation space $\Rep(Q,\mathbf{v})$ consists of all representations of the given quiver $Q$ with dimension vector $\mathbf{v}$.   
Consider the Jordan quiver, and let $\mathbf{v}=n$. 
Then $\Rep(Q,n) 
= \spec(\mathbb{C}[V_r]) \cong \mathbb{C}^{n^2}$.  
If $\mathbf{v} = (n,m)$ for the $A_2$-quiver,
then $\Rep(Q,(n,m)) = \spec(\mathbb{C}[V_a]) \cong \mathbb{C}^{mn}$.   
 
Studying isomorphism classes of representations of a quiver    
with a prescribed dimension vector amounts to a change-of-basis of the vector space at each vertex, i.e.,   
it amounts to a certain natural quotient of the representation space by a group, 
and the study of the orbit structure is equivalent to equivariant geometry 
in geometric invariant theory (GIT; see \S\ref{subsection:background-quivers-GIT}).   

%\subsubsection{Filtered quiver representations}
%\label{subsubsection:filtered-quiver-rep}
%{\color{red} here } 

\subsubsection{Filtered ADHM data}
\label{subsubsection:filtered-ADHM-data}
In this section, we generalize the constructions in \cite{MR1604167}, giving a refined analogue of quiver representations called filtered quiver representations.
A \textbf{\color{Brown}filtered quiver representation} is a finite quiver with a fixed filtration of vector spaces at each vertex 
 whose linear map associated to each arrow (of the quiver) preserves the filtration. 
 Furthermore, there is a natural group action by 
   the parabolic group $\mathfrak{P}_{\mathbf{v}}$ of $\mathfrak{G}_{\mathbf{v}}=\prod_{\iota} \GL_{v_{\iota}}$ acting on and preserving filtered vector spaces.
    
Let $Q=(Q_0,Q_1)$, where cycles are allowed, and let 
$Q_1^{\op}$ be the same set of arrows as $Q_1$ but in opposite orientation. 
Let $\overline{Q}=(Q_0,Q_1\sqcup Q_1^{\op})$, a \textbf{\color{Brown} double quiver}.   
Write $\overline{Q}_1 := Q_1\sqcup Q_1^{\op}$.   
Let $\mathbf{v}=(v_{\iota}) \in \mathbb{N}^{Q_0}$  and let 
$V=(V_{\iota})_{\iota\in Q_0}$   
be a collection of vector spaces such that 
$\dim V_{\iota} =v_{\iota}$ for each $\iota \in Q_0$.   
For 
$V^1=(V_{\iota}^1)_{\iota\in Q_0}$  
and 
$V^2 =(V_{\iota}^2)_{\iota\in Q_0}$, 
define  
\[
L(V^1, V^2) := \displaystyle{ \bigoplus_{\iota\in Q_0}}  \Hom(V_{\iota}^1, V_{\iota}^2)
\quad 
\mbox{ and }
\quad 
E(V^1,V^2) := \displaystyle{ \bigoplus_{a\in \overline{Q}_1 }  } \Hom(V_{ta}^1, V_{ha}^2).  
\]  
Given 
$B=(B_a)_{a\in \overline{Q}_1 } \in E(V^1,V^2)$ 
and  
$C=(C_a)_{a\in \overline{Q}_1 } \in E(V^2,V^3)$,  
we also define 
\[ 
CB := \left(\displaystyle{ \sum_{ta=\iota} C_a B_{a^{\op}} } \right)_{\iota\in Q_0}\in L(V^1,V^3),
\]  
where  
$a^{\op}$ has the same endpoints as $a$ but is in reverse orientation.

Now,  choose a sequence $\gamma^1,\gamma^2,\ldots, \gamma^N =\mathbf{v}\in\mathbb{N}^{Q_0}$ of dimension vectors such that for all $k$ and  $\iota\in Q_0$, $\gamma_{\iota}^k\leq \gamma_{\iota}^{k+1}$. 
For each $\iota\in Q_0$, we get a filtration of $\mathbb{C}^{v_{\iota}}$: 
\[ 
\{ 0\} \subseteq \mathbb{C}^{\gamma_{\iota}^1}\subseteq \mathbb{C}^{\gamma_{\iota}^2}\subseteq \ldots \subseteq \mathbb{C}^{v_{\iota}},  
\] 
where each $\mathbb{C}^l$ is spanned by $l$ basis elements of $\mathbb{C}^{v_{\iota}}$. 
Define 
\[ 
\text{F}^{\bullet}\Rep(Q,\mathbf{v},\mathbf{w}) := 
\text{F}^{\bullet}E(V,V)\oplus L(W,V)\oplus L(V,W), 
\]  
where 
\[ 
\text{F}^{\bullet}E(V,V) 
:= 
\displaystyle{ \bigoplus_{a\in \overline{Q}_1 } 
\text{F}^{\bullet}\Hom(V_{ta},V_{ha})}\] 
 such that if 
$B_a\in \text{F}^{\bullet}\Hom(V_{ta},V_{ha})$ where $a\in Q_1$, then $B_a$ preserves the fixed sequence of vector spaces at every level, i.e.,  
$B_a(\mathbb{C}^{\gamma_{ta}^k} )\subseteq \mathbb{C}^{\gamma_{ha}^k}$ for every $k$, 
and if 
$B_c\in \text{F}^{\bullet}\Hom(V_{tc},V_{hc})$, where $c=a^{\op} \in Q_1^{\op}$, then 
$\text{F}^{\bullet}\Hom(V_{tc},V_{hc})$ is the dual of $\text{F}^{\bullet}\Hom(V_{ta},V_{ha})$ such that the trace map 
\[ 
\text{F}^{\bullet}\Hom(V_{ta},V_{ha}) \times \text{F}^{\bullet}\Hom(V_{tc},V_{hc})\stackrel{\tr}{\longrightarrow}\mathbb{C}, \hspace{4mm}
(B_a,B_{c})\mapsto \Tr(B_aB_{c})  
\] 
is a nondegenerate pairing.   
We denote 
\[ 
C= (A,B)\in \text{F}^{\bullet}E(V,V), 
\:\:
\mbox{ where }
\:\:
A\in \bigoplus_{a\in Q_1}\text{F}^{\bullet}\Hom(V_{ta},V_{ha})
\:\:
\mbox{ and }
\:\: 
B\in 
\bigoplus_{c\in Q_1^{\op}}\text{F}^{\bullet}\Hom(V_{tc},V_{hc}), 
\]  
$i\in L(W,V)$, and $j\in L(V,W)$. 
We call an element of 
$\text{F}^{\bullet}\Rep(Q,\mathbf{v},\mathbf{w})$ a 
{\bf\color{Brown} filtered ADHM datum},  
while the filtered representation space  
$\text{F}^{\bullet}\Rep(Q,\mathbf{v},\mathbf{w})$ is called 
\textbf{\color{Brown}filtered ADHM data}.

\subsection{Invariant and semi-invariant polynomials}
\label{subsection:inv-semi-inv}

Invariant and semi-invariant polynomials play a fundamental role in classical and geometric invariant theory (\S\ref{subsection:background-quivers-GIT}). 
In fact, studying orbit spaces precisely amounts to describing invariant and semi-invariant polynomials.

\subsection{Moment maps and complete intersection}\label{subsection:intro-moment-maps-complete-intersection}  
Moment maps arise in symplectic geometry as a tool to construct conserved quantities.  
The action of a Lie group $G$ on a vector space $X$ is induced to the cotangent bundle $T^*X$ of $X$. 
Taking the derivative of the group action induces an infinitesimal action 
$\mathfrak{g}\stackrel{a}{\rightarrow} TX$ on $X$, given  
by tangent vectors. By dualizing this action, one obtains the moment map $T^*X \stackrel{\mu}{\rightarrow} \mathfrak{g}^*$, where $\mu = a^*$. 
If the zero fiber of $\mu$ is a complete intersection, then this means $\mu^{-1}(0)$ has the expected number of irreducible components, with $\mu^{-1}(0)$ having dimension $2\dim X-\dim \mathfrak{g}$. 
% 
% 
%has finitely many connected components, the ring of $G$-invariant polynomials is a finitely generated algebra, 
%and studying the geometry over the zero locus is equivalent to studying over the $\xi$-fiber, for any $\xi\in \mathfrak{g}^*$. 
% \cite{MR1458060}  

\subsubsection{Moment maps for filtered ADHM data}
\label{subsubsection:moment-map-filt-ADHM}
Recall filtered ADHM data from \S\ref{subsubsection:filtered-ADHM-data}. 

Let $\epsilon: \overline{Q}_1  \rightarrow \{\pm 1 \}$, where 
\[   
\epsilon(a)=\begin{cases}   
 \:\:\: 1 & \mbox{ if }a\in Q_1,   \\  
     -1 & \mbox{ if }a\in Q_1^{\op}. \\   
\end{cases}  
\] 
Then $\epsilon C \in \text{F}^{\bullet}E(V^1,V^2)$ is defined as 
$(\epsilon C)_a = \epsilon(a)C_a$ for $a\in \overline{Q}_1$. 
We define a symplectic form $\omega$ on 
$\text{F}^{\bullet}\Rep(Q,\mathbf{v},\mathbf{w})$ by 
\[ 
\omega((C,i,j),(C',i',j'))=\Tr(\epsilon C C')+\Tr(ij'-i'j), 
\quad 
\mbox{ where } 
\quad  
\Tr(A)=\sum_k \Tr(A_k). 
\]  
The product $\mathfrak{P}_{\mathbf{v}}:= \prod_{\iota\in Q_0}P_{v_{\iota}}$ of parabolic groups acts on 
$\text{F}^{\bullet}\Rep(Q,\mathbf{v},\mathbf{w})$ 
via 
\[ 
p\circ (C,i,j)=(p C p^{-1},pi,jp^{-1}), 
\]  
which preserves the symplectic form $\omega$ on the filtered ADHM data. 
The $\mathfrak{P}_{\mathbf{v}}$-action induces the moment map 
\begin{equation}
\label{eq:moment-map-filt-ADHM}
\mu_{\mathfrak{P}_{\mathbf{v}}}: 
T^*\text{F}^{\bullet}\Rep(Q,\mathbf{v},\mathbf{w}) \rightarrow 
\lie(\mathfrak{P}_{\mathbf{v}})^*=\bigoplus_{\iota\in Q_0}\mathfrak{p}_{v_{\iota}}^*, 
\quad 
\mbox{ where } 
\quad 
(C,i,j)\mapsto \epsilon CC+ij  
\end{equation}
and $\mathfrak{p}_{v_{\iota}}=\Lie(P_{v_{\iota}})$.

\subsection{Geometric invariant theory}\label{subsection:background-quivers-GIT}  
The orbit space for a $G$-action on a variety $X$ may not exist since $X/G$ may not necessarily be an algebraic variety. 
In order to remedy this, 
invariant polynomials are used to construct new and interesting varieties called affine quotients; 
%For example,    
%consider the action of the complex torus $\mathbb{C}^*$ acting on the complex plane $\mathbb{C}^2$ by $t.(x,y)=(tx,ty)$.  
%All the $\mathbb{C}^*$-orbits are $1$ (complex) dimensional, except when $(x,y)=(0,0)$; the origin in the complex plane is a fixed %point. 
%The natural quotient map    
%$\mathbb{C}^2 \rightarrow \mathbb{C}^2/\mathbb{C}^*$ does not exist since $\mathbb{C}^2/\mathbb{C}^*$ is not an orbit space.  
%In order to remedy this, we define what is known as the affine space whose algebra is generated by invariant polynomials, i.e.,   
%polynomials that remain invariant under the action of the group.      
%In the above example, $\mathbb{C}^2 /\!\!/ \mathbb{C}^* = \spec(\mathbb{C})=\mathbb{C}^0$   
%since there are no invariant polynomials.   
one may think of the affine quotient as being generated by orbit closures (or alternatively, closed orbits).   
Procesi in~\cite{MR0419491} proved that invariant polynomials of quiver varieties  
arise as traces of oriented paths in characteristic zero, 
and Donkin in~\cite{MR1185589} and \cite{MR1259609} later proved an analogous result in characteristic $p$.

To produce other interesting (and sometimes projective) varieties,   
one uses GIT techniques (cf.~\cite{MR1304906,MR2537067}), 
where one aims to produce all semi-invariant polynomials for various character of the group.   
%In the example involving the $1$-dimensional torus, consider the character $\chi:t\mapsto t$. 
%Let $f(x,y)=x$. 
%Then since $f(t.(x,y))=tx$, $f$ is a $\chi$-semi-invariant polynomial. 
%Similarly, for $g(x,y)=y$, we see that $g$ is a $\chi$-semi-invariant polynomial. Thus 
%we obtain     
%$\mathbb{C}^2/\!\!/_{\chi}\mathbb{C}^*   
%= \Proj(\mathbb{C} \oplus \mathbb{C}\cdot x\oplus \mathbb{C}\cdot y \oplus \ldots )
%=\Proj(\mathbb{C}[x,y])=\mathbb{P}_{\mathbb{C}^0}^{1}$, the projective line. 
Derksen--Weyman in ~\cite{MR1758750}, Schofield--van den Bergh in ~\cite{MR1908144},  
and Domokos--Zubkov in ~\cite{MR1825166} used long exact sequences of the Ringel resolution, %~\cite{MR0422350}, 
representation-theoretic techniques, and combinatorial techniques, respectively, to give a strategy to produce semi-invariant polynomials for quiver representations.

 Recall that given a $\mathfrak{G}$-action on $\mathfrak{X}$, 
 affine and GIT quotients are 
 \[ 
 \mathfrak{X}/\!\!/\mathfrak{G} := \spec(\mathbb{C}[\mathfrak{X}]^{\mathfrak{G}}) 
 \quad 
 \mbox{ and }
 \quad 
 \mathfrak{X}/\!\!/_{\chi}\mathfrak{G} := \Proj\left(\bigoplus_{i=0}^{\infty} \mathbb{C}[\mathfrak{X}]^{\mathfrak{G},\chi^i}\right), 
 \] 
respectively, where $\chi:\mathfrak{G}\rightarrow \mathbb{C}^*$ is a character of $\mathfrak{G}$.  
 For more detail, see \cite{MR1304906,MR2537067}.

%
%
%
%
% 
% IMPORTANT:  \subsection{Hamiltonian reduction and Gr\"obner basis}\label{subsection:background-Ham-reduction-gb}  
% IMPORTANT:  Hamiltonian reduction appears from the $G_{\beta}$-equivariant moment map construction by taking the preimage of $0$ as seen in Section~\ref{subsection:intro-moment-maps-complete-intersection} and quotienting out the locus by $G_{\beta}$. Gr\"obner bases are thoroughly discussed in \cite{MR1322960} and \cite{MR1287608},  which is essentially a procedure using a notion of $s$-polynomials to make the generating set of an ideal ``bigger" (at least in the setting of Hamiltonian reduction).   
% IMPORTANT:  I study $\mu^{-1}(0)$ via its ideal in $\mathbb{C}[T^*X]$.  
%  
% 

\subsubsection{Hamiltonian reduction of filtered ADHM data}
\label{subsubsection:quiver-varieties}  
Recall $\mu_{\mathfrak{P}_{\mathbf{v}}}$ defined in \eqref{eq:moment-map-filt-ADHM}. 
The locus $\mu_{\mathfrak{P}_{\mathbf{v}}}=0$ is called \textbf{\color{Brown}filtered ADHM equation}. We say 
\[ 
\mathfrak{M}_0^{\text{F}^{\bullet}} = \mathfrak{M}_0^{\text{F}^{\bullet}}(\mathbf{v},\mathbf{w}):=\mu_{\mathfrak{P}_{\mathbf{v}}}^{-1}(0)/\!\!/\mathfrak{P}_{\mathbf{v}}
=\spec\left( \mathbb{C}[\mu_{\mathfrak{P}_{\mathbf{v}}}^{-1}(0) ]^{\mathfrak{P}_{\mathbf{v}} }\right) 
\] 
 is the \textbf{\color{Brown} filtered affine quotient}, 
and    
\[ 
\mathfrak{M}_{\chi}^{\text{F}^{\bullet}} = 
\mathfrak{M}_{\chi}^{\text{F}^{\bullet}}(\mathbf{v},\mathbf{w}) := 
\mu_{\mathfrak{P}_{\mathbf{v}}}^{-1}(0)^s/\mathfrak{P}_{\mathbf{v}} = 
\Proj\left(\displaystyle{\bigoplus_{i\geq 0}} \mathbb{C}[\mu_{\mathfrak{P}_{\mathbf{v}}}^{-1}(0) ]^{\mathfrak{P}_{\mathbf{v}},\chi^i} \right) 
\]  
is called a \textbf{\color{Brown}filtered quiver variety}.

In the case when $Q$ is the Jordan quiver 
\begin{equation}
\label{eqn:framed-Jordan-quiver}
  \xymatrix{   
  \stackrel{1}{\bullet} \ar@(rd,ru)_{r},  
	       }   
\end{equation}
then the double quiver $\overline{Q}$ is 
\begin{equation}
\label{eqn:GS-quiver}
  \xymatrix{   
  \stackrel{1.}{\bullet} \ar@(rd,ru)_{r} \ar@(lu,ld)_{r^{\op}},  
	       }  
\end{equation}
Let $\mathbf{v}=n$ and $\mathbf{w}=1$. 
We impose the complete standard filtration of vector spaces on the representation $V_1=\mathbb{C}^n$ at vertex $1$ to obtain that  $\mu_B^{-1}(0)/B\cong T^*(\widetilde{\mathfrak{g}}\times \mathbb{C}^n/\GL_n(\mathbb{C}))$ (cf.~\cite[Prop.~3.2]{Nevins-GSresolutions}, \cite[Prop.~1.1]{Im-rss-locus-GS-resolution}).

%Let $G$ be a reductive group over $\mathbb{C}$ and let $\mathfrak{g}=\Lie(G)$. 
%Let $P$ and $B$ be a parabolic and Borel subgroup of $G$, respectively, $\mathfrak{p}=\Lie(P)$ and $\mathfrak{b}=\Lie(B)$.   

\section{Invariants and semi-invariants of filtered quiver representations}\label{section:results-filtered-semi-invariants}    
We obtain the following in \cite[Thm.~5.1.2]{Im-doctoral-thesis}: 
\begin{theorem}[Im]
\label{theorem:filtered-quiver-vars-finite-Dynkin-type}
Consider a quiver $Q=(Q_0,Q_1)$ of finite Dynkin type with dimension vector $\mathbf{v} = (n,\ldots, n)$. Impose the complete standard filtration of vector spaces on the representation at each vertex of the quiver.  
Let $\mathfrak{U}$  be the product of maximal unipotent subgroups of the product $\mathfrak{B}=B^{Q_0}$ of standard Borels. 
Then 
$\mathbb{C}[\mathfrak{b}^{Q_1}]^{\mathfrak{U}}\cong\mathbb{C}[\mathfrak{t}^{Q_1}]$.  
\end{theorem}
    
Theorem~\ref{theorem:filtered-quiver-vars-finite-Dynkin-type} says that invariant polynomials cannot arise from off-diagonal coordinates of $\mathfrak{b}$. 
Instead, only the eigenvalues of $\mathfrak{b}^{Q_1}$ 
remain invariant under the $\mathfrak{U}$-action.

Next, we state \cite[Thm.~5.2.12]{Im-doctoral-thesis}.
\begin{theorem}[Im]
\label{theorem:filtered-quiver-vars-affine-Dynkin-type} 
Consider an affine quiver $\widetilde{A}_r$ with a framing, and let  
$\mathbf{v} = (n,\ldots, n,m)$, the dimension vector. 
Assume the complete standard filtration of vector spaces on the representation at each vertex of $\widetilde{A}_r$ (except at the framed vertex). 
Let $\mathfrak{U}$  be the product of maximal unipotent subgroups of the product $\mathfrak{B}=B^{r}$ of standard Borels. 
Then the invariant subring 
$\mathbb{C}[\mathfrak{b}^{\oplus r+1}\oplus M_{n\times m}]^{\mathfrak{U}}$ 
has finitely-many generators.  
\end{theorem} 

We give an explicit description of the subalgebra in 
Theorem~\ref{theorem:filtered-quiver-vars-affine-Dynkin-type}, 
thus generalizing  
\cite[Thm.~13.3]{MR1489234}, which states that for the maximal unipotent group 
$U$ of the standard Borel $B \subseteq \GL_n(\mathbb{C})$, 
the algebra $\mathbb{C}[M_{n\times m}]^U$ 
is generated by bideterminants of standard Young bitableaux $(D|E)$, 
where each row of $D$ has the form $p, p+1, \ldots, n$ for a suitable $p$, $1\leq p\leq n$. 
Note that one may associate a bideterminant to a Young bitableau in a natural way, i.e., it is a product of minors, each of which is determined by $i$-th row of $D$ and $i$-th row of $E$. 
 
Also see \cite[Thm.~1.1, Thm.~1.3]{Im-two-pathways} as further extensions of Theorem~\ref{theorem:filtered-quiver-vars-finite-Dynkin-type} and Theorem~\ref{theorem:filtered-quiver-vars-affine-Dynkin-type}.

\section{Hamiltonian reduction of the Borel moment map}
\label{section:Ham-red-Borel}  
Throughout this section, we will restrict to the Jordan quiver.   
Let $\mathbf{v}=n$ and $\mathbf{w}=1$. Let $\text{F}^{\bullet}$ be the complete standard filtration of vector spaces. 
Then $\text{F}^{\bullet}\Rep(Q,n,1)=T^*(\mathfrak{b}\times \mathbb{C}^n)$, which is associated to the moment map 
\[ 
\mu_B:T^*(\mathfrak{b}\times \mathbb{C}^n)\rightarrow \mathfrak{b}^*, 
\quad 
\mbox{ where } 
\quad 
(r,s,i,j)\mapsto [r,s]+ij. 
\] 
Let $\mu_B^{-1}(0)^{\rss}$ be the locus of points whose eigenvalues of $r$ are pairwise distinct. 
Let $\Delta_n= \{(r_1,\ldots, r_n,0,\ldots, 0): r_{\iota} =r_{\gamma} 
\mbox{ for some }\iota \not=\gamma\}$. Then \cite[Thm.~1.5]{Im-rss-locus-GS-resolution} says that: 
\begin{theorem}[Im]
\label{thm:P-surjection}
The map 
\[ 
P:\mu_B^{-1}(0)^{\rss}\twoheadrightarrow\mathbb{C}^{2n}\setminus \Delta_n
\quad  
\mbox{ given by }
\quad 
(r,s,i,j)\mapsto (r_{11},\ldots, r_{nn},s_{11}',\ldots, s_{nn}')
\]  
is a regular, well-defined surjection separating $B$-orbit closures. 
\end{theorem}
Furthermore,  we have  \cite[Thm.~1.6]{Im-rss-locus-GS-resolution}: 
\begin{theorem}[Im] 
\label{thm:P-induces-isom}
The surjection $P$ in Theorem~\ref{thm:P-surjection} descends to an isomorphism 
\[ 
\mu_B^{-1}(0)^{\rss}/\!\!/B \cong \mathbb{C}^{2n}\setminus \Delta_n
\] 
of algebraic varieties. 
\end{theorem} 
We prove that the components of $\mu_B$ form a regular sequence using a certain monomial ordering and the resulting initial ideal.  

In the process of proving Theorem~\ref{thm:P-surjection} and Theorem~\ref{thm:P-induces-isom}, 
rational $B$-invariant polynomials appearing as traces of products of matrices were discovered (cf.~\cite[\S6.2]{Im-doctoral-thesis}): for $1\leq \iota\leq n$  and $1\leq \gamma< \nu \leq n$,  
\begin{equation}
\label{eqn:B-inv-functions}
\begin{split}
F_{\iota}(r,s,i,j) &=  \tr\left(j L^{\iota} i  \right), 
			\\
G_{\iota}(r,s,i,j)&=
  \tr \left( L^{\iota} s\right),  
			\\ 
H_{\iota}(r,s,i,j)&=\tr(L^{\iota}r),  
%= e_{\iota}^* \: r\: e_{\iota}  
			\\
K_{\gamma\nu}(r,s,i,j) &= [\tr((L^{\nu}-L^{\gamma})r)]^{-1}.  
%   = (e_{\nu}^* (r-r_{\gamma\gamma}\I)e_{\nu})^{-1} 
\end{split} 
\end{equation}
 
These rational functions should play an important role in the construction of the affine quotient $\mu_B^{-1}(0)/\!\!/B$ (see \S\ref{subsection:affine-GIT-open} for more detail).

\section{Hamiltonian reduction of the parabolic moment map}
\label{section:Ham-red-parabolic} 
In this section, we continue to work with the Jordan quiver.   
Let $\mathbf{v}=n$ and $\mathbf{w}=1$. 
Here, rather than working with a complete flag as in \S\ref{section:Ham-red-Borel}, we work with a partial standard flag 
\begin{equation}
\label{eqn:flag-max-5-vs}
\text{F}^{\bullet}: \mathbb{C}^{\alpha_1}\subseteq \mathbb{C}^{\alpha_1+\alpha_2}\subseteq \ldots \subseteq \mathbb{C}^{\alpha_1+\ldots +\alpha_{\ell}}, 
\quad 
\mbox{ where }
\quad  
\ell \leq 5. 
\end{equation}
Then we have 
\cite[Thm.~1.1]{Im-Scrimshaw-parabolic}:  
 \begin{theorem}[Im--Scrimshaw]
 \label{theorem:Grothendieck--Springer-resolution-complete-intersection} 
Let $\alpha = (\alpha_1, \dotsc, \alpha_{\ell})$ such that $\ell \leq 5$.
Let $P$ be the parabolic subgroup of $\GL_n(\CC)$ with block size vector $\alpha$.
The components of $\mu_P$ form a complete intersection.  
\end{theorem}
We also explicitly describe the irreducible components in \cite[Thm.~1.2]{Im-Scrimshaw-parabolic}.

\section{Open problems}
\label{section:results-problems}   

We will now list some research problems related to filtered ADHM data. 

\subsection{(Semi-)invariant polynomials for a general quiver}
As a result of the $B$-invariant functions in \eqref{eqn:B-inv-functions}, let us now focus on the 
$n$ rational functions $L^{\iota}$ enumerated by $\iota = 1, \ldots, n$. 
For $l_k(r)=r-r_{kk}\I$,  
\[  
M^{\iota}:=  
 \left(\tr\left(\prod_{1\leq k\leq n, k\not= \iota} l_k(r)\right)\right)  L^{\iota} 
= \prod_{1\leq k\leq n, k\not= \iota} l_k(r), 
\] 
which is an operator whose matrix entries are zero except in coordinates 
$(\mu,\nu)$ for $\mu\leq \iota$ and $\nu\geq \iota$.  
One should think of these operators acting on elements in $\mathfrak{b}$ as killing off columns of a matrix, 
or as creating new coordinate entries from matrices in $\mathfrak{b}^*=\mathfrak{gl}_n/\mathfrak{u^+}$  
such that powers of the trace of the product of these matrices give the desired invariant polynomials for the filtered affine quotient setting.

%The path algebra $\mathbb{C}Q$ of a quiver is the $\mathbb{C}$-vector space whose basis is the set of all paths in $Q$. 
%In the first example, $\mathbb{C}Q=\mathbb{C}[r]$; it is an infinite-dimensional vector space with generators $r^k$ where $k\geq 0$,  
%with $r^0 := e_1$ is defined to be the trivial path at vertex $1$.  
%The path algebra of the second quiver is $\mathbb{C}Q=\mathbb{C}e_1 \oplus \mathbb{C}a \oplus \mathbb{C}e_2$, a three dimensional %vector space.  
%We define the product $pq$ of two paths $p$ and $q$     
%to be the concatenation of paths if the tail of $p$ equals the head of $q$, and zero otherwise.  

%The notion of Jordan canonical matrices does not quite make sense in the study of $B$-orbits on $\mathfrak{b}$ 
%since the $B$-action fixes the diagonal entries of any   
%$r \in \mathfrak{b}$ and $r$ cannot be put into a Jordan canonical form using only the matrices in $B$.  
%(although an analogous canonical form does exist for $B$-action on $\mathfrak{b}$).   

%     I will introduce a new technique through a notion called filtered quivers and give strategy on how one may study parabolic equivariant geometry on filtered vector spaces. 

% There are at least three notions of quotienting: 
% $\mathfrak{gl}_n/G$, $\mathfrak{gl}_n/\!\!/G$, $\mathfrak{gl}_n/\!\!/_{\chi}G$, in the hopes of constructing orbit space, affine quotient, and geometric invariant theory (GIT) quotient, respectively.    
			  
\begin{problem}
\label{problem:inv-Kronecker-Jordan-quiver}
Describe $\mathbb{C}[\mathfrak{b}^{\oplus k}]^{U_n\times U_n}$ for $k\geq 3$, and $\mathbb{C}[\mathfrak{b}^{\oplus \ell}]^{U_n}$ for $\ell\geq 2$. 
\end{problem} 
  
By making progress on 
Problem~\ref{problem:inv-Kronecker-Jordan-quiver}, together with Theorem~\ref{theorem:filtered-quiver-vars-affine-Dynkin-type}, 
one can then describe the unipotent invariant subring for a filtered ADHM data for \textit{any} quiver $Q$.

\subsection{Affine and GIT quotients}
\label{subsection:affine-GIT-open}

Let $Q$ be the Jordan quiver. 
By clearing the denominators of the rational functions in \eqref{eqn:B-inv-functions}, they are \textit{some} of the generators of 
$\mathbb{C}[\mu_B^{-1}(0)/\!\!/B]$. More generally, we have:      
\begin{problem}\label{problem:relations-among-generators}    
Find generators and relations for 
$\mu_{\mathfrak{P}_{\mathbf{v}}}^{-1}(0)/\!\!/\mathfrak{P}_{\mathbf{v}}$ and 
$\mu_{\mathfrak{P}_{\mathbf{v}}}^{-1}(0)/\!\!/_{\chi}\mathfrak{P}_{\mathbf{v}}$
to describe interesting quotients.  
\end{problem}    

Problem~\ref{problem:relations-among-generators} 
is of interest for constructing new moduli spaces, 
and it remains an open problem to study the algebro-geometric structure of these quotients for various quivers, filtrations, and stability conditions, i.e., for different characters. 

\subsection{Birational morphisms}
\label{subsection:birational-mor} 

\begin{problem}\label{problem:variational-GIT}   
Construct birational morphisms between $\mu_{\mathfrak{P}_{\mathbf{v}}}^{-1}(0)/\!\!/_{\chi}\mathfrak{P}_{\mathbf{v}}$ 
and 
$\mu_{\mathfrak{P}_{\mathbf{v}}}^{-1}(0)/\!\!/_{\chi'}\mathfrak{P}_{\mathbf{v}}$ for two characters $\chi$ and $\chi'$ of $\mathfrak{P}_{\mathbf{v}}$. 
\end{problem}  

Problem~\ref{problem:variational-GIT} is known as variational GIT.  
Fixing a character of the group and constructing GIT quotients is known as a polarization on the variety, and changing one polarization to another is known as wall-crossing. 
Since a parabolic group has many more characters than $\GL_n(\mathbb{C})$, this problem now becomes even more interesting than a reductive group setting.  
Thus what happens when we fix a character of the parabolic group and the phenomenon that arises when we cross a wall? 
Are the resulting GIT quotients birational or isomorphic?    

Even when we restrict the quiver to the Jordan quiver, studying variational GIT and constructing Fourier automorphism, thus giving an isomorphism between certain open loci of GIT quotients, are doable yet a difficult task 
(one needs to choose various $1$-parameter subgroups and show that certain points are unstable with respect to a fixed character).

\subsection{Complete intersection for $\mu_P^{-1}(0)$}
\label{subsection:Borel-moment-map-ci-open}
It remains to prove that when $\ell >5$ for $\ell$ in \eqref{eqn:flag-max-5-vs}, the associated parabolic moment map is flat, i.e.,  
    $\mu_P^{-1}(0)$ is a complete intersection or equivalently, the components of $\mu_P$ form a regular sequence. 
    Furthermore, one should use the well-known fact that  
    $\mu_P$ is flat if and only if 
		$\dim \mu_P^{-1}(0)/\!\!/P = 2(\dim (\mathfrak{p}\times \mathbb{C}^n)-\dim \mathfrak{p})$ 
if and only if  
\[ \codim \{ y\in \mathfrak{p}\times \mathbb{C}^n: \dim G_y = k \} \geq k \hspace{2mm}\mbox{ for all }k\geq 1.  
\]

Results in \cite{Im-rss-locus-GS-resolution} coincide with the classical notion that the trace of an oriented cycle 
of a quiver, as well as the trace of a path that begin and end at a framed vertex, is an invariant function (cf.~\cite{MR1834739,MR958897,MR1623674}). 
These strategies are appliable to Nakajima's affine and quiver varieties.  
This leads us to believe, together with \cite[Proof of Prop.~8.2.1]{MR2210660}, that $\mathbb{C}[\mu_P^{-1}(0)]^P$ is generated by traces of products of matrices.

\begin{problem}
The filtered ADHM equation $\mu_P$ is a complete intersection for $\ell >5$. 
\end{problem}

\subsection{Complete intersection for the general quiver}
\label{subsection:complete-intersection-filtered}
%the $\chi_i$-locus to $\chi_{n+1-i}^{-1}$-locus.     
Regular sequences for the filtered ADHM equation play an important role in affine and GIT quotients.  
A set of polynomials $f_1,\ldots, f_N$ is   
$\mathbb{C}[x_1,\ldots, x_m]$-\textbf{\color{Brown}regular} if the scheme defined by the vanishing locus of $f_1,\ldots, f_N$ form a complete intersection.   
That is, the variety has the expected dimension $m-N$, and in such a case, 
we say that the components form a \textbf{\color{Brown}regular sequence}.      
The affine quotient of a variety is generated by invariant polynomials, 
and GIT quotients are generated by semi-invariant polynomials.   
In the case that the components of a moment map form a regular sequence, the affine and the GIT quotients are finitely generated; this means we only need to produce finitely-many invariant and semi-invariant polynomials since the quotients will have the expected dimension.

\begin{problem}\label{problem:components-regular-sequence}    
Describe when components of $\mu_{\mathfrak{P}_{\mathbf{v}}}$ form a regular sequence.   
\end{problem}

\subsection{Comparing filtered ADHM data and quivers with relations}\label{subsection:results-quivers-with-relations}      
There is a close relation to quivers with relations but examples show that they are not equivalent. It would be interesting 
to check under what conditions are irreducible components of quivers with relations normal. 
In the case that components of quivers with relations are normal, the function theory for quivers with relations can be extended to all of the component, and, thus, is comparable to filtered ADHM data.

%A relation of a quiver $Q$ is a subspace of $\mathbb{C}Q$ spanned by linear combinations of paths having a common source and a common target,   
%with each combination of paths having length at least $2$. 
%Quivers with relations are thought of as imposing relations on the composition of arrows, which is equivalent to quotienting the path algebra by the appropriate ideal.  
 
%It is a classical fact that any commutative and noncommutative algebra can be represented using quivers with relations. 
%That is, it can be thought of as path algebra with relations.   
%In fact, the category of representations of quivers is equivalent to the category of finite dimensional $\mathbb{C}Q$-modules.    
%In the category of representations of quivers with relations, all arrows represent linear maps which must satisfy the imposed relations. 

%%%%%%%%%%%%%%%%%%%%%%%%%%%%%%%%%%%%%%%%%%%%%%%%%%%%%%%%%%%%%%%%%%%%%%%%
%
%
%
%
%%%%%%%%%%%%%%%%%%%%%%%%%%%%%%%%%%%%%%%%%%%%%%%%%%%%%%%%%%%%%%%%%%%%%%%%
  
\subsection{Hilbert schemes}
\label{subsection:Hilbert-schemes}  
 Writing 
 $T^*(\mathfrak{b}\times \mathbb{C}^n)$ as $\mathfrak{b}\times \mathfrak{b}^*\times \mathbb{C}^n \times (\mathbb{C}^n)^*$, 
 it contains the set of all quadruples of the form $(r,s,i,j)$, 
 where $s$ takes the form of lower triangular matrices (technically, $\mathfrak{b}^*=\mathfrak{gl}_n/\mathfrak{u}$ where $\mathfrak{u}$ is the set of nilpotent matrices in $\mathfrak{b}$),   
 $i$ is a vector, and $j$ is a covector.  
 %its associated filtered quiver representation is: 
% \[      
% \xymatrix@R-1pc{ 
% 							& & \mathbb{C}^n     	\ar@{-}[d]   \\     
% 							& & \vdots 						\ar@{-}[d]   \\      
% 							& & \mathbb{C}^1     	\ar@{-}[d]   \\    
% \mathbb{C}^1 	& & \mathbb{C}^0     								\\     
% \stackrel{1}{\circ} \ar@/^/[rr]| i && \stackrel{2}{\bullet} \ar@/^/[ll]|{j}  \ar@(r,d)[]^{s}  \ar@(r,u)[]_{r}
% }      
% \]   
%    where the framing is assigned at vertex $1$.   
    Framing means that there is no group action on the vector space at that vertex.  
    Restricting to those points satisfying the relation $[r,s]+ij=0$ should remind the experts of the Hilbert scheme 
    $(\mathbb{C}^2)^{[n]}$ of $n$ points on the complex plane. 
    
    \begin{problem}
    Construct 
    $\mu_P^{-1}(0)/\!\!/_{\chi}P\stackrel{}{\twoheadrightarrow} \mu_P^{-1}(0)/\!\!/P$ and relate it to the Hilbert--Chow morphism 
    $(\mathbb{C}^2)^{[n]}\stackrel{}{\twoheadrightarrow} \mathbb{C}^{2n}/S_n$, where $S_n$ is a permutation group of $n$ letters.  
\end{problem}    
     
% 
%Although I was able to prove that such inequality holds for all general matrices 
%up to $\mathfrak{b}_4$ by hand (i.e., these are matrices with indeterminates as coordinates), 
%  the strategy has been difficult to generalize.     
%  My second strategy has been to use deformation theory techniques and a certain lemma in \cite{MR1735775} 
%  by linking $\mu_B^{-1}(0)$ to another variety that is a complete intersection, but it seems that I need to use Gr\"obner basis if I want to finish in this way.  
%    Finally, I use Gr\"obner degeneration and Gr\"obner basis strategy by imposing certain weights on the variables of the ambient space (note that we have $2(n(n+1)/2+n)$ variables), apply Gr\"obner degeneration and when looking at the initial ideal, 
%    the problem becomes somewhat simplified. 
%    Before imposing any form of degeneration techniques, Macaulay2 (M2) could not finish computing that $\mu_B^{-1}(0)\subseteq T^*(\mathfrak{b}_4\times \mathbb{C}^4)$ is a regular sequence, but after imposing the degeneration,  
%    M2 was able to show that $\mu_B^{-1}(0)\subseteq T^*(\mathfrak{b}_5\times \mathbb{C}^5)$ is a complete intersection within two minutes. Gr\"obner techniques seem possible but still very difficult 
%    for larger matrices.  
%    To prove the general case, 
%    I need to write down a number of syzygy matrices whose entries consist of polynomials or zero.
%    It is possible that more experimental data is needed to understand the patterns. 

\subsection{Idempotents in filtered quiver representations and quivers with relations}\label{subsection:idempotents-filtered-quiv-quiv-reln}   
The $L^{\iota}$'s in \cite{Im-rss-locus-GS-resolution} (also see \cite{Im-doctoral-thesis}) are $n$ rational idempotents. They are pairwise orthonormal and their sum equals the identity matrix.   
These idempotents are all upper triangular matrices whose coordinate are rational functions.    
In \S\ref{subsection:results-quivers-with-relations}, 
we compare filtered ADHM data with quivers with relations. 
  In view of the Borel subalgebra corresponding to the Jordan quiver embedded in a quiver with relations, whose 
 quiver has $n$ vertices and each vertex has a trivial path (also known as idempotents in its path algebra), it is highly likely that these idempotent matrices are directly related to the $n$ idempotents appearing in quivers with relations.   
 
\subsection{Modified rational Cherednik algebras} 
\label{subsection:rational-Cherednik-alg}

In \cite{MR2210660}, the authors consider the spherical subalgebra of $\mathfrak{gl}_n$-type and realize it as a quantum Hamiltonian reduction of the algebra $\mathcal{D}(\mathfrak{g})$ of polynomial differential operators on $\mathfrak{g}$. Since studying the Hamiltonian reduction of $\mathcal{D}(\mathfrak{g})$ with respect to $P$ is the same as performing the Hamiltonian reduction of $\mathcal{D}(\mathfrak{g}\times \mathbb{P}^n)$ with respect to $G$ (acting diagonally on $\mathfrak{g}\times \mathbb{P}^n$), this leads us to an analogous investigation for parabolic groups: 

\begin{problem}
\label{problem:cherednik-algebras-D-modules}
Construct a parabolic analog of the Cherednik algebra and realize it as a quantization of $\mathcal{D}(\mathfrak{p})$. 
\end{problem}

\subsection{Parabolic $\mathcal{D}$-modules} 
\label{subsection:parabolic-D-modules}
In a similar spirit as in \cite{MR2444305}, 
we believe that the quantization of the Hilbert scheme associated to $\mu_P^{-1}(0)/\!\!/_{\chi}P$ may be realized as a microlocalization of a modified rational Cherednik algebra. 
Furthermore, Gan--Ginzburg construct a Lagrangian subvariety $\mathfrak{M}_{\text{nil}}$ in $\mathfrak{M}$ in the classical case of $GL_n(\mathbb{C})$, and then studies a category of holonomic $\mathcal{D}$-modules whose characteristic variety is contained in $\mathfrak{M}_{\text{nil}}$.

\begin{problem}\label{problem:category-characteristic-variety}
Construct a category of $\mathcal{D}$-modules whose characteristic variety is contained in the parabolic analogue $\mathfrak{M}_{\text{nil}}^P$ and compare its simple objects to Lusztig's (parabolic) character sheaves (cf.~\cite{MR2074987,MR2124170}). 
\end{problem}

\subsection{Derived category of coherent sheaves on $\mu_P^{-1}(0)/\!\!/_{\chi}P$}
\label{subsection:derived-filtered}

There is an equivalence between derived categories of coherent sheaves on $(\mathbb{C}^2)^{[n]}$ and finitely-generated modules over a noncommutative crepant resolution of $\mathbb{C}^{2n}/S_n$ (cf.~\cite{MR1918676}). 

\begin{problem}
Construct an equivalence between derived categories of coherent sheaves on $\mu_P^{-1}(0)/\!\!/_{\chi}P$ and the category of finitely-generated modules over a modified rational Cherednik algebra. 
\end{problem}

%Next, we make a distinction between quiver varieties and quiver representations.   
%When one refers to quiver varieties, one usually means Nakajima quiver varieties, i.e., the Hamiltonian reduction of a double quiver variety twisted by a nontrivial character.   
%In this paper, we will only consider quiver representations.  

\subsection{Generalized Grothendieck--Springer resolution}\label{subsection:generalized-grothendieck-springer-resolution}

 It is a classical result in representation theory that the Grothendieck--Springer resolution $\widetilde{\mathfrak{g}}\twoheadrightarrow \mathfrak{g}$, where $(x,\mathfrak{b})\mapsto x$, is generically finite and dominant of degree $|W|$, where $W$ is the Weyl group. 
 
 \begin{problem}
 \label{problem:extend-GS-resoln-setting}
 Extend $\mu_{\mathfrak{P}_{\mathbf{v}}}$ from the Jordan quiver to a general quiver as it is interesting to study the resolution from a filtered quiver variety point of view.
 \end{problem}

Directly studying the parabolic orbits on the Grothendieck--Springer resolution is very difficult, and surprisingly, no one seems to have studied such geometric objects from the quiver variety point of view.

Filtered ADHM data represent a new technique for studying objects in geometric representation theory and algebraic geometry.  
They are nice because of their concrete nature and because of many applications in mathematics. 
Using filtered quiver variety techniques, 
we hope that new and different insights to 
important varieties will be identified and discovered.

\appendix
\bibliography{inv-and-semi-inv-of-all-filtered-quivers}   

\def\cprime{$'$} \def\cprime{$'$} \def\cprime{$'$} \def\cprime{$'$}
\begin{thebibliography}{10}

\bibitem{MR618892}
Walter Borho and Robert MacPherson.
\newblock Repr\'{e}sentations des groupes de {W}eyl et homologie d'intersection
  pour les vari\'{e}t\'{e}s nilpotentes.
\newblock {\em C. R. Acad. Sci. Paris S\'{e}r. I Math.}, 292(15):707--710,
  1981.

\bibitem{MR2838836}
Neil Chriss and Victor Ginzburg.
\newblock {\em Representation theory and complex geometry}.
\newblock Modern Birkh\"auser Classics. Birkh\"auser Boston Inc., Boston, MA,
  2010.
\newblock Reprint of the 1997 edition.

\bibitem{MR2772068}
Alastair Craw.
\newblock Quiver flag varieties and multigraded linear series.
\newblock {\em Duke Math. J.}, 156(3):469--500, 2011.

\bibitem{Crawley-Boevey-rep-quivers}
William Crawley-Boevey.
\newblock Lectures on representations of quivers.
\newblock \url{http://www.maths.leeds.ac.uk/~pmtwc/quivlecs.pdf}, 1992.

\bibitem{MR1834739}
William Crawley-Boevey.
\newblock Geometry of the moment map for representations of quivers.
\newblock {\em Compositio Math.}, 126(3):257--293, 2001.

\bibitem{MR1758750}
Harm Derksen and Jerzy Weyman.
\newblock Semi-invariants of quivers and saturation for
  {L}ittlewood-{R}ichardson coefficients.
\newblock {\em J. Amer. Math. Soc.}, 13(3):467--479 (electronic), 2000.

\bibitem{MR1825166}
Matyas Domokos and Alexandr~N. Zubkov.
\newblock Semi-invariants of quivers as determinants.
\newblock {\em Transform. Groups}, 6(1):9--24, 2001.

\bibitem{MR1185589}
Stephen Donkin.
\newblock Invariants of several matrices.
\newblock {\em Invent. Math.}, 110(2):389--401, 1992.

\bibitem{MR1259609}
Stephen Donkin.
\newblock Polynomial invariants of representations of quivers.
\newblock {\em Comment. Math. Helv.}, 69(1):137--141, 1994.

\bibitem{MR2210660}
Wee~Liang Gan and Victor Ginzburg.
\newblock Almost-commuting variety, $\mathscr{D}$-modules, and {C}herednik
  algebras.
\newblock {\em IMRP Int. Math. Res. Pap.}, pages 26439, 1--54, 2006.
\newblock With an appendix by Ginzburg.

\bibitem{MR1649626}
Victor Ginzburg.
\newblock Geometric methods in the representation theory of {H}ecke algebras
  and quantum groups.
\newblock In {\em Representation theories and algebraic geometry ({M}ontreal,
  {PQ}, 1997)}, volume 514 of {\em NATO Adv. Sci. Inst. Ser. C Math. Phys.
  Sci.}, pages 127--183. Kluwer Acad. Publ., Dordrecht, 1998.
\newblock Notes by Vladimir Baranovsky [V. Yu. Baranovski{\u\i}].

\bibitem{Ginzburg-Nakajima-quivers}
Victor Ginzburg.
\newblock Lectures on {N}akajima's quiver varieties.
\newblock \url{http://arxiv.org/pdf/0905.0686v2}, 2009.

\bibitem{MR1489234}
Frank~D. Grosshans.
\newblock {\em Algebraic homogeneous spaces and invariant theory}, volume 1673
  of {\em Lecture Notes in Mathematics}.
\newblock Springer-Verlag, Berlin, 1997.

\bibitem{MR1918676}
Mark Haiman.
\newblock Vanishing theorems and character formulas for the {H}ilbert scheme of
  points in the plane.
\newblock {\em Invent. Math.}, 149(2):371--407, 2002.

\bibitem{Im-two-pathways}
Mee~Seong Im.
\newblock Semi-invariants of filtered quiver representations with at most two
  pathways.
\newblock arXiv:1409.0702.

\bibitem{Im-doctoral-thesis}
Mee~Seong Im.
\newblock {\em On semi-invariants of filtered representations of quivers and
  the cotangent bundle of the enhanced {G}rothendieck-{S}pringer resolution}.
\newblock ProQuest LLC, Ann Arbor, MI, 2014.
\newblock Thesis (Ph.D.)--University of Illinois at Urbana-Champaign.

\bibitem{Im-rss-locus-GS-resolution}
Mee~Seong Im.
\newblock The regular semisimple locus of the affine quotient of the cotangent
  bundle of the {G}rothendieck-{S}pringer resolution.
\newblock {\em J. Geom. Phys.}, 132:84--98, 2018.

\bibitem{Im-Scrimshaw-parabolic}
Mee~Seong Im and Travis Scrimshaw.
\newblock The regularity of almost-commuting partial {G}rothendieck--{S}pringer
  resolutions and parabolic analogs of {C}alogero--{M}oser varieties.
\newblock arXiv:1812.02283.

\bibitem{MR2444305}
Masaki Kashiwara and Rapha\"{e}l Rouquier.
\newblock Microlocalization of rational {C}herednik algebras.
\newblock {\em Duke Math. J.}, 144(3):525--573, 2008.

\bibitem{MR2525917}
Mikhail Khovanov and Aaron~D. Lauda.
\newblock A diagrammatic approach to categorification of quantum groups. {I}.
\newblock {\em Represent. Theory}, 13:309--347, 2009.

\bibitem{MR2763732}
Mikhail Khovanov and Aaron~D. Lauda.
\newblock A diagrammatic approach to categorification of quantum groups {II}.
\newblock {\em Trans. Amer. Math. Soc.}, 363(5):2685--2700, 2011.

\bibitem{MR958897}
Lieven Le~Bruyn and Claudio Procesi.
\newblock Semisimple representations of quivers.
\newblock {\em Trans. Amer. Math. Soc.}, 317(2):585--598, 1990.

\bibitem{MR2124170}
G.~Lusztig.
\newblock Parabolic character sheaves. {II}.
\newblock {\em Mosc. Math. J.}, 4(4):869--896, 981, 2004.

\bibitem{MR1035415}
George Lusztig.
\newblock Canonical bases arising from quantized enveloping algebras.
\newblock {\em J. Amer. Math. Soc.}, 3(2):447--498, 1990.

\bibitem{MR1182165}
George Lusztig.
\newblock Canonical bases arising from quantized enveloping algebras. {II}.
\newblock {\em Progr. Theoret. Phys. Suppl.}, (102):175--201 (1991), 1990.
\newblock Common trends in mathematics and quantum field theories (Kyoto,
  1990).

\bibitem{MR1088333}
George Lusztig.
\newblock Quivers, perverse sheaves, and quantized enveloping algebras.
\newblock {\em J. Amer. Math. Soc.}, 4(2):365--421, 1991.

\bibitem{MR1623674}
George Lusztig.
\newblock On quiver varieties.
\newblock {\em Adv. Math.}, 136(1):141--182, 1998.

\bibitem{MR1775358}
George Lusztig.
\newblock Quiver varieties and {W}eyl group actions.
\newblock {\em Ann. Inst. Fourier (Grenoble)}, 50(2):461--489, 2000.

\bibitem{MR1758244}
George Lusztig.
\newblock Semicanonical bases arising from enveloping algebras.
\newblock {\em Adv. Math.}, 151(2):129--139, 2000.

\bibitem{MR2074987}
George Lusztig.
\newblock Parabolic character sheaves. {I}.
\newblock {\em Mosc. Math. J.}, 4(1):153--179, 311, 2004.

\bibitem{MR1304906}
David Mumford, John Fogarty, and Frances Kirwan.
\newblock {\em Geometric invariant theory}, volume~34 of {\em Ergebnisse der
  Mathematik und ihrer Grenzgebiete (2) [Results in Mathematics and Related
  Areas (2)]}.
\newblock Springer-Verlag, Berlin, third edition, 1994.

\bibitem{MR1302318}
Hiraku Nakajima.
\newblock Instantons on {ALE} spaces, quiver varieties, and {K}ac-{M}oody
  algebras.
\newblock {\em Duke Math. J.}, 76(2):365--416, 1994.

\bibitem{MR1604167}
Hiraku Nakajima.
\newblock Quiver varieties and {K}ac-{M}oody algebras.
\newblock {\em Duke Math. J.}, 91(3):515--560, 1998.

\bibitem{MR1711344}
Hiraku Nakajima.
\newblock {\em Lectures on {H}ilbert schemes of points on surfaces}, volume~18
  of {\em University Lecture Series}.
\newblock American Mathematical Society, Providence, RI, 1999.

\bibitem{Nevins-GSresolutions}
Thomas Nevins.
\newblock Stability and {H}amiltonian reduction for {G}rothendieck-{S}pringer
  resolutions.
\newblock
  \url{http://www.math.uiuc.edu/~nevins/papers/b-hamiltonian-reduction-2011-0316.pdf},
  2011.

\bibitem{MR2537067}
Peter~E. Newstead.
\newblock Geometric invariant theory.
\newblock In {\em Moduli spaces and vector bundles}, volume 359 of {\em London
  Math. Soc. Lecture Note Ser.}, pages 99--127. Cambridge Univ. Press,
  Cambridge, 2009.

\bibitem{MR0419491}
Claudio Procesi.
\newblock The invariant theory of {$n\times n$} matrices.
\newblock {\em Advances in Math.}, 19(3):306--381, 1976.

\bibitem{Rouquier-2-Kac-Moody-algebras}
Raphael Rouquier.
\newblock 2-{K}ac-{M}oody algebras.
\newblock \url{http://arxiv.org/pdf/0812.5023v1}, 2008.

\bibitem{MR1908144}
Aidan Schofield and Michel van~den Bergh.
\newblock Semi-invariants of quivers for arbitrary dimension vectors.
\newblock {\em Indag. Math. (N.S.)}, 12(1):125--138, 2001.

\end{thebibliography}
\bibliographystyle{plain}

\end{document}